\documentclass[reqno]{amsproc}
\usepackage{bbm}
\usepackage{amssymb,hyperref,color}

\newtheorem{theorem}{Theorem}[section]

\theoremstyle{remark}

\numberwithin{equation}{section}
\allowdisplaybreaks

\title[]
{Morrey spaces over the unit circle cannot be renormed\\ to become rearrangement-invariant}

\author[O. Karlovych]{Oleksiy Karlovych}
\address{
Centro de Matem\'atica e Aplica\c{c}\~oes\\
Departamento de Matem\'atica\\
Faculdade de Ci\^encias e Tecnologia\\
Universidade Nova de Lisboa\\
Quinta da Torre\\
2829--516 Caparica\\
Portugal} \email{oyk@fct.unl.pt}
\author[E. Shargorodsky]{Eugene Shargorodsky}
\address{%
Department of Mathematics\\
King's College London\\
Strand, London WC2R 2LS\\
United Kingdom}
\email{eugene.shargorodsky@kcl.ac.uk}
\begin{document}
\begin{abstract}
Let $1\le p<\infty$ and $0<\lambda<1$. We consider the classical Morrey space
$L^{p,\lambda}(\mathbb{T})$ over the unit circle $\mathbb{T}$. We show that 
there are equimeasurable functions $f,g:\mathbb{T}\to\mathbb{R}$ such that 
$g\in L^{p,\lambda}(\mathbb{T})$ but $f\notin L^{p,\lambda}(\mathbb{T})$.
This implies that the the space $L^{p,\lambda}(\mathbb{T})$ cannot be 
renormed to become rearrangement-invariant.
\end{abstract}
\keywords{Morrey space, rearrangement-invariant space.}
\subjclass{Primary 46E30}
\maketitle
\section{Introduction and the main result}
Let $\mathbb{T}$ be the unit circle and $m$ be the Lebesgue measure on 
$\mathbb{T}$ normalised so that $m(\mathbb{T})=1$. For $1\le p<\infty$ and 
$0\le\lambda<1$, the Morrey space $L^{p,\lambda}(\mathbb{T})$ is the collection  
of all measurable functions $f:\mathbb{T}\to\mathbb{C}$ such that
\[
\|f\|_{L^{p,\lambda}}:=\sup_{\omega\subset\mathbb{T}}
\left(
\frac{1}{m(\omega)^\lambda}\int_\omega |f(\zeta)|^p dm(\zeta)
\right)^{1/p}<\infty,
\]
where the supremum is taken over all arcs $\omega\subset\mathbb{T}$. 
It is well known that $L^{p,\lambda}(\mathbb{T})$ is a Banach function
space in the sense of \cite[Ch.~1, Definition~1.1]{BS88}
(cf. \cite[Sections~1.2.1 and 1.2.4]{SDH20}, where Morrey spaces
over $\mathbb{R}^n$ are considered).
 
Recall that two measurable functions $f,g:\mathbb{T}\to\mathbb{C}$ are
said to be equimeasurable if
\[
m\{\zeta\in\mathbb{T}:|f(\zeta)|>\lambda\}
=
m\{\zeta\in\mathbb{T}:|g(\zeta)|>\lambda\}
\]
for all $\lambda\ge 0$.
A Banach function space $X(\mathbb{T})$ is said to be rearrangement-invariant
if for every pair of equimeasurable functions $f,g:\mathbb{T}\to\mathbb{C}$ one
has $\|f\|_X=\|g\|_X$. We are interested in the question whether 
the Morrey spaces $L^{p,\lambda}(\mathbb{T})$ are rearrangement-invariant.

One can replace the unit circle $\mathbb{T}$ by the real line $\mathbb{R}$
and define the Morrey space $L^{p,\lambda}(\mathbb{R})$ in a similar fashion.
The behaviour of functions in the Morrey space $L^{p,\lambda}(\mathbb{R})$
near infinity may be quite complicated. In particular, it follows from 
\cite[Lemma~6.1]{HS17} that $L^{p,\lambda}(\mathbb{R})$ is not embedded into 
$L^1(\mathbb{R})+L^\infty(\mathbb{R})$ if $0<\lambda<1$. Therefore, in view of 
\cite[Ch.~1, Theorem~6.6]{BS88}, the Morrey space 
$L^{p,\lambda}(\mathbb{R})$ is not rearrangement-invariant if $0<\lambda<1$. 

The above complications do not occur in the case of the unit circle,
and one has to use a different argument to show that
$L^{p,\lambda}(\mathbb{T})$ is not
rearrangement-invariant even after renorming. The following is the main result of the paper.
\begin{theorem}
Let $1\le p<\infty$ and $0<\lambda<1$. Then there exist equimeasurable functions
$f,g:\mathbb{T}\to\mathbb{R}$ such that $g\in L^{p,\lambda}(\mathbb{T})$ but
$f\notin L^{p,\lambda}(\mathbb{T})$.
\end{theorem}
This theorem immediately implies that the Morrey space 
$L^{p,\lambda}(\mathbb{T})$ cannot be renormed so that it becomes
rearrangement-invariant. This result is probably well known to experts, but we have not
been able to find an appropriate reference. The main aim of the paper is to provide such a reference. 
\section{Proof}
Let $0 < \varepsilon < \min\left\{\lambda/2, 1 - \lambda\right\}$. 
For $-\pi \le \theta \le \pi$, consider the functions
\begin{align*}
f\left(e^{i\theta}\right) &:= 
\left\{\begin{array}{lll}
   n^{(1 - \lambda + \varepsilon)/p},   
   & \frac{1}{n + 1} < \theta < \frac1n , & n = 16,17,\dots, 
   \\
    0  & \text{otherwise},
\end{array} \right.
\\
g\left(e^{i\theta}\right) &:= 
\left\{\begin{array}{lll}
   n^{(1 - \lambda + \varepsilon)/p} ,  
   & \frac1{\sqrt{n}} - \frac{1}{n(n + 1)} < \theta < \frac1{\sqrt{n}} , 
   & n =16,17,\dots, \\
    0  & \text{otherwise}.
\end{array}\right.
\end{align*}
Clearly, $f$ and $g$ are equimeasurable. 

Since 
\[
\frac{1}{n+1}<\theta<\frac{1}{n}\le\frac{1}{16},
\quad
n=16,17,\dots,
\]
we have 
\[
n>\frac{1}{\theta}-1\ge\frac{1}{2\theta},
\quad
n=16,17,\dots
\]
and
\[
f\left(e^{i\theta}\right) \ge (2\theta)^{-(1 - \lambda + \varepsilon)/p}, 
\quad
0 < \theta <\frac{1}{16}.
\]
Therefore, for every $t\in(0,1/16)$ and 
$\omega_t := \left\{e^{i\theta}: 0 < \theta < t\right\}$, we have
\begin{align*}
\frac{1}{m(\omega_t)^\lambda}\int_{\omega_t} |f(\zeta)|^p dm(\zeta) 
& \ge  
\frac{1}{(t/(2\pi))^\lambda} \frac{2^{-(1-\lambda+\varepsilon)}}{2\pi}
\int_0^t \theta^{-(1 - \lambda + \varepsilon)} d\theta
\\
&= 
\frac{(4\pi)^{\lambda - 1}}{2^\varepsilon} 
t^{-\lambda}\, 
\frac{1}{\lambda - \varepsilon}\, t^{\lambda - \varepsilon} 
= 
\frac{(4\pi)^{\lambda - 1}}{2^\varepsilon(\lambda - \varepsilon)}\, 
t^{-\varepsilon} \to +\infty
\end{align*}
as $t \to 0+$. Hence $f \not\in L^{p,\lambda}(\mathbb{T})$. 

On the other hand, we are going to show that $g \in L^{p,\lambda}(\mathbb{T})$. 
Let $\mathbb{T}_+ := \left\{e^{i\theta}:0 < \theta < \frac{1}{4}\right\}$. 
Since $g(\zeta) = 0$ for $\zeta \in \mathbb{T}\setminus\mathbb{T}_+$, one has
\[
\sup_{\omega\subset\mathbb{T}}
\left(\frac{1}{m(\omega)^\lambda}\int_\omega |g(\zeta)|^p dm(\zeta)\right)^{1/p} 
=
\sup_{\omega\subset\mathbb{T}_+}
\left(\frac{1}{m(\omega)^\lambda}\int_\omega |g(\zeta)|^p dm(\zeta)\right)^{1/p} .
\]
Take any $\omega\subset\mathbb{T}_+$. We only need to consider arcs $\omega$ that 
intersect at least one of the arcs 
\[
\gamma_n := \left\{e^{i\theta} :  
\frac1{\sqrt{n}} - \frac{1}{n(n + 1)} < \theta < \frac1{\sqrt{n}}\right\}, 
\quad
n= 16, 17,\dots. 
\]
Set
\begin{align*}
n_1 & := \max\left\{
n \in \mathbb{N}: 
\inf_{e^{it} \in \omega} t < \frac{1}{\sqrt{n}}
\right\} , 
\\
n_0 & := \min\left\{
n \in \mathbb{N}: 
\sup_{e^{it} \in \omega} t > \frac{1}{\sqrt{n}} - \frac{1}{n(n + 1)}
\right\} .
\end{align*}
Since $n_1$ corresponds to the leftmost arc $\gamma_n$ intersecting $\omega$ 
while $n_0$ corresponds to the rightmost such arc, $n_0 \le n_1$.

If $n_0 = n_1$, then
\begin{align*}
\frac{1}{m(\omega)^\lambda}\int_\omega |g(\zeta)|^p dm(\zeta) 
& = 
\frac{1}{m(\omega)^\lambda}
\int_{\omega\cap\gamma_{n_0}} |g(\zeta)|^p dm(\zeta)  
\\ 
&\le 
\frac{1}{m(\omega\cap\gamma_{n_0})^\lambda}
\int_{\omega\cap\gamma_{n_0}} |g(\zeta)|^p dm(\zeta) 
\\
&= 
m(\omega\cap\gamma_{n_0})^{1 - \lambda} n_0^{1 - \lambda + \varepsilon} 
\le 
m(\gamma_{n_0})^{1 - \lambda} n_0^{1 - \lambda + \varepsilon} 
\\
& = 
\big(2\pi n_0(n_0 + 1)\big)^{\lambda - 1} n_0^{1 - \lambda + \varepsilon} 
< 
n_0^{\lambda - 1 + \varepsilon} < 1.
\end{align*}

Suppose now $n_0 < n_1$. Then
\begin{align*}
m(\omega) 
&\ge 
\frac{1}{\sqrt{n_0}} - \frac{1}{n_0(n_0 + 1)} - \frac{1}{\sqrt{n_1}} 
\\
& \ge 
\frac12 \left(\frac{1}{\sqrt{n_0}} 
- 
\frac{1}{\sqrt{n_1}}\right) 
+ 
\frac12 \left(\frac{1}{\sqrt{n_0}} 
- 
\frac{1}{\sqrt{n_0 + 1}}\right) 
- 
\frac{1}{n_0(n_0 + 1)} 
\\
& = 
\frac12 \left(
\frac{1}{\sqrt{n_0}} - \frac{1}{\sqrt{n_1}}
\right) 
+ 
\frac{1}{2\sqrt{n_0}\sqrt{n_0 + 1}(\sqrt{n_0 + 1} + \sqrt{n_0})} 
- 
\frac{1}{n_0(n_0 + 1)} 
\\
& >  
\frac12 \left(\frac{1}{\sqrt{n_0}} - \frac{1}{\sqrt{n_1}}\right) 
+
\frac{1}{4(n_0 + 1)\sqrt{n_0}} - \frac{1}{n_0(n_0 + 1)} 
\ge 
\frac12 \left(\frac{1}{\sqrt{n_0}} - \frac{1}{\sqrt{n_1}}\right),
\end{align*}
where the last inequality holds because $n_0 \ge 16$. Further,
\begin{align*}
2\pi \int_\omega |g(\zeta)|^p dm(\zeta)  
& \le 
\sum_{n = n_0}^{n_1} n^{1 - \lambda + \varepsilon} \frac{1}{n(n + 1)} 
< 
\sum_{n = n_0}^{n_1} n^{-1 - \lambda + \varepsilon} 
\\
& < 
4  \sum_{n = n_0}^{n_1} \int_n^{n + 1} x^{-1 - \lambda + \varepsilon} dx 
= 
4 \int_{n_0}^{n_1 + 1} x^{-1 - \lambda + \varepsilon} dx 
\\
& < 
8 \int_{n_0}^{n_1} x^{-1 - \lambda + \varepsilon} dx 
< 
8 \int_{n_0}^{n_1} x^{-1 - \lambda/2} dx 
\\
& = 
\frac{16}{\lambda}\left(n_0^{-\lambda/2} - n_1^{-\lambda/2}\right) .
\end{align*}
Hence
\begin{align*}
\frac{1}{m(\omega)^\lambda}\int_\omega |g(\zeta)|^p dm(\zeta) 
& < 
2^\lambda \left(\frac{1}{\sqrt{n_0}} - \frac{1}{\sqrt{n_1}}\right)^{-\lambda}
\frac{16}{2\pi \lambda}\left(n_0^{-\lambda/2} - n_1^{-\lambda/2}\right) 
\\
& = 
\frac{2^{\lambda + 3}}{\pi \lambda} 
\left(\left(\frac{n_1}{n_0}\right)^{1/2} - 1\right)^{-\lambda} 
\left(\left(\frac{n_1}{n_0}\right)^{\lambda/2} - 1\right) .
\end{align*}
Let 
\[
\phi_\lambda(y) := 
\left(y^{1/2} - 1\right)^{-\lambda} \left(y^{\lambda/2} - 1\right) , 
\quad y > 1 .
\]
Since $\phi_\lambda$ is continuous on $(1, \infty)$ and
\[
\lim_{y \to 1} \phi_\lambda(y) = 0 , 
\qquad 
\lim_{y \to \infty} \phi_\lambda(y) = 1 ,
\]
one has
\[
M_\lambda := \sup_{y \in (1, \infty)} \phi_\lambda(y) < \infty ,
\]
and
\[
\frac{1}{m(\omega)^\lambda}
\int_\omega |g(\zeta)|^p dm(\zeta) 
< 
\frac{2^{\lambda + 3} M_\lambda}{\pi \lambda}\, .
\]
So, $g \in L^{p,\lambda}(\mathbb{T})$.
\qed

The above method allows one to easily construct equimeasurable functions $f,g:\mathbb{R}\to\mathbb{R}$ such that 
$g\in L^{p,\lambda}(\mathbb{R})$ but $f\notin L^{p,\lambda}(\mathbb{R})$.

\subsection*{Funding}
This work is funded by national funds through the FCT - Funda\c{c}\~ao para a 
Ci\^encia e a Tecnologia, I.P., under the scope of the projects UIDB/00297/2020 
(\url{https://doi.org/10.54499/UIDB/00297/2020})
and UIDP/ 00297/2020 
(\url{https://doi.org/10.54499/UIDP/00297/2020})
(Center for Mathematics and Applications).
\bibliographystyle{abbrv}
\bibliography{OKES23}
\end{document}